\documentclass[11pt,leqno,twoside]{article}

\usepackage{amsfonts,amsmath,amsthm,amssymb}
\usepackage{enumerate}
\usepackage{graphics,graphicx,subfigure}

\usepackage{color}
\usepackage{todonotes}
\usepackage{cancel}
\usepackage{url}
\usepackage{hyperref}
\usepackage{makeidx}
\usepackage{showidx}
\usepackage{multicol}        
\usepackage{xspace}
\usepackage{stmaryrd}        
\usepackage{pifont}          
\usepackage{fancybox}        
\usepackage{bm}

 \usepackage{fancyhdr}
\theoremstyle{plain}
 \theoremstyle{definition}
 \newtheorem{lem}{Lemma}
 \newtheorem{defn}[lem]{Definition}
 \newtheorem{thm}[lem]{Theorem}
 \newtheorem{prop}[lem]{Proposition}
 \newtheorem{cor}[lem]{Corollary}
 \newtheorem{notn}[lem]{Notations}
 \newtheorem{pb}[lem]{Problem}
 \newtheorem{form}[lem]{Formulae}
 
 \newtheorem*{rk}{Remark}
 \newtheorem*{com}{Comment}
 \newtheorem*{ex}{Example}
 \theoremstyle{remark}

 \newcommand{\blem}{\begin{lem}}
 \newcommand{\elem}{\end{lem}}
 \newcommand{\bdefn}{\begin{defn}}
 \newcommand{\edefn}{\end{defn}}
 \newcommand{\bthm}{\begin{thm} }
 \newcommand{\ethm}{\end{thm}}
 \newcommand{\bprop}{\begin{prop}}
 \newcommand{\eprop}{\end{prop}}
 \newcommand{\bcor}{\begin{cor}}
 \newcommand{\ecor}{\end{cor}}
 \newcommand{\bnotn}{\begin{notn}}
 \newcommand{\enotn}{\end{notn}}
 \newcommand{\bpb}{\begin{pb}}
 \newcommand{\epb}{\end{pb}}
 \newcommand{\bform}{\begin{form}}
 \newcommand{\eform}{\end{form}}
 \newcommand{\brk}{\begin{rk}}
 \newcommand{\erk}{\end{rk}}
 \newcommand{\bcom}{\begin{com}}
 \newcommand{\ecom}{\end{com}}
 \newcommand{\bex}{\begin{ex}}
 \newcommand{\eex}{\end{ex}}
 \newcommand{\bpf}{\begin{proof}}
 \newcommand{\epf}{\end{proof}}

\newcommand{\cC}{\mathcal{C}}

\newcommand{\cE}{\mathcal{E}}

\newcommand{\cK}{\mathcal{K}}

\newcommand{\cM}{\mathcal{M}}

\newcommand{\cP}{\mathcal{P}}

\newcommand{\bN}{\mathbb{N}}

\newcommand{\bR}{\mathbb{R}}

\newcommand{\be}{\begin{equation}}
\newcommand{\ee}{\end{equation}}
\newcommand{\bal}{\begin{align}}
\newcommand{\eal}{\end{align}}
\newcommand{\ba}{\begin{align*}}
\newcommand{\ea}{\end{align*}}
\newcommand{\bmx}{\begin{matrix}}
\newcommand{\emx}{\end{matrix}}
\newcommand{\bbmx}{\begin{bmatrix}}
\newcommand{\ebmx}{\end{bmatrix}}
\newcommand{\bpmx}{\begin{pmatrix}}
\newcommand{\epmx}{\end{pmatrix}}
\newcommand{\bvmx}{\begin{vmatrix}}
\newcommand{\evmx}{\end{vmatrix}}

\newcommand{\ul}{\underline}

\newcommand{\wh}{\widehat}
\newcommand{\wt}{\widetilde}
\newcommand{\f}{\frac}

\newcommand{\inc}{\subseteq}

\newcommand{\setm}{\setminus}
\newcommand{\tto}{\longrightarrow}

\newcommand{\tr}{\mathrm{tr}}

\newcommand{\minimize}[1]{\underset{#1}{\rm minimize}\,}

\newcommand{\la}{\lambda}

\newcommand{\eps}{\varepsilon}
  
\title{\vspace{-20mm}
Instances of Computational Optimal Recovery:\\
Refined Approximability Models
\medskip\hrule height 1.2pt \vspace{-6mm}}

\author{Simon Foucart\footnote{S. F. is partially supported by NSF grants DMS-1622134 and DMS-1664803, and also acknowledges the NSF grant CCF-1934904.} --- Texas A\&M University}
\date{\vspace{-6mm}\rule{100mm}{0.8pt}}

\newcommand\shorttitle{Instances of Computational Optimal Recovery:
Refined Approximability Models}
\newcommand\authors{S. Foucart}

\begin{document}
\maketitle

\vspace{-15mm}
\begin{abstract}
Models based on approximation capabilities have recently been studied in the context of Optimal Recovery.
These models, however,
are not compatible with overparametrization,
since model- and data-consistent functions could then be unbounded.
This drawback motivates the introduction of refined approximability models featuring an added boundedness condition.
Thus, two new models are proposed in this article:
one where the boundedness applies to the target functions (first type)
and one where the boundedness applies to the approximants (second type).
For both types of model,
optimal maps for the recovery of linear functionals are first described on an abstract level
before their efficient constructions are addressed.
By exploiting techniques from semidefinite programming,
these constructions are explicitly carried out on a common example involving polynomial subspaces of $\cC[-1,1]$.
\end{abstract}



\vspace{-5mm}
\begin{center}
\rule{100mm}{0.8pt}
\end{center}

\section{Introduction}\vspace{-2mm}

The objective of this article is to uncover practical methods for the optimal recovery of functions available through observational data when the underlying models based on approximability  allow for overparametrization.
To clarify this objective and its various challenges,
we start with some background on traditional Optimal Recovery.
Typically, an unknown function $f$ defined on a domain~$D$ is observed through point evaluations $y_1 = f(x_1),\ldots,y_m=f(x_m)$
at distinct points $x_1,\ldots,x_m \in D$.
More generally, an unknown object $f$,
simply considered as an element of a normed space $X$,
is observed through
\be
\label{Obs}
y_i = \ell_i(f),
\qquad i \in [1:m],
\ee
where $\ell_1,\ldots,\ell_m$ are linear functionals defined on $X$.
We assume here that these data are perfectly accurate 
 --- we refer to the companion article \cite{EttFou} for the incorporation of observation error.
The data is summarized as $y = L(f)$, where the linear map $L: X \to \bR^m$ is called observation operator.
Based on the knowledge of $y \in \bR^m$,
the task is then to recover a quantity of interest $Q(f)$,
where throughout this article $Q: X \to \bR$ is assumed to be a linear functional.
The recovery procedure can be viewed as a map $R$ from $\bR^m$ to $\bR$,
with no concern for its practicability at this point.

Besides the observational data (which is also called {\sl a posteriori} information),
there is some {\sl a priori} information coming from an educated belief about the properties of realistic $f$'s.
It translates into the assumption that $f$ belongs to a model set $\cK \inc X$.
The choice of this model set is of course critical.
When the $f$'s indeed represent functions,
it is traditionally taken as the unit ball with respect to some norm that characterizes smoothness.
More recently,
motivated by parametric partial differential equations,
a model based on approximation capabilities has been proposed in~\cite{BCDDPW}.
Namely, given a linear subspace $V$ of $X$ and a threshold $\eps > 0$, 
it is defined as
\be
\label{ApproxMod}
\cK = \cK_{V,\eps} :=\{
f \in X: {\rm dist}_X(f,V) \le \eps \}.
\ee
This model set is also implicit in many numerical procedures
and in machine learning.

Whatever the selected model set,
the performance of the recovery procedure $R: \bR^m \to \bR$ is measured in a worst-case setting via the (global) error of $R$ over $\cK$, i.e.,
\be
e_{\cK,Q}(L,R) := \sup_{f \in \cK} |Q(f) - R(L(f))|.
\ee
Obviously, one is interested in optimal recovery maps $R^{\rm opt}: \bR^m \to \bR$ minimizing
this worst-case error, i.e., such that
\be
e_{\cK,Q}(L,R^{\rm opt}) = \inf_{R: \bR^m \to \bR} e_{\cK,Q}(L,R).
\ee
This infimum is called the intrinsic error of the observation map $L$ (for $Q$ over $\cK$).
It is known, at least since Smolyak's doctoral dissertation \cite{Smo},
that there is a linear functional among the optimal recovery maps as soon as
the set $\cK$ is symmetric and convex,
see e.g. \cite[Theorem 4.7]{TractI} for a proof.
The practicality of such a linear optimal recovery map is not automatic, though.
For the approximability set~\eqref{ApproxMod},
Theorem 3.1 of \cite{DFPW} revealed that such a linear optimal recovery map takes the form $R^{\rm opt}: y \in \bR^m \mapsto \sum_{i=1}^m a_i^{\rm opt} y_i$,
where $a^{\rm opt} \in \bR^m$ is a solution to
\be
\label{ResultDFPW}
\minimize{a \in \bR^m} \bigg\| Q - \sum_{i=1}^m a_i \ell_i \bigg\|_{X^*}
\qquad \mbox{subject to  }
\sum_{i=1}^m a_i \ell_i(v) = Q(v)
\quad \mbox{for all }v \in V,
\ee
an optimization problem that can be solved for $X = \cC(D)$
in exact form when the observation functionals are point evaluations (see \cite{DFPW})
and in approximate form when they are arbitrary linear functionals (see \cite{EttFou} or Subsection \ref{SecCompI} below).

The approximability set \eqref{ApproxMod}, however, presents some important restrictions.
Suppose indeed that there is some nonzero $v \in \ker(L) \cap V$.
Then, 
for a given $f_0 \in \cK$ observed through $y = L(f_0) \in \bR^m$,
any $f_t := f_0 + t v$, $t \in \bR$, is both model-consistent (i.e., $f_t \in \cK$) and data-consistent (i.e., $L(f_t) = y$),
so that the local error at $y$ of any recovery map $R: \bR^m \to \bR$
satisfies
\be
e^{\rm loc}_{\cK,Q}(L,R(y)) := \sup_{\substack{f \in \cK \\ L(f) = y}}
|Q(f) - R(y)| \ge
\sup_{t \in \bR} |Q(f_t)-R(y)|
= \sup_{t \in \bR} |(Q(f_0)-R(y)) + t Q(v)|,
\ee
which is generically infinite.
Thus, for the optimal recovery problem to make sense under the approximability model \eqref{ApproxMod},
one must assume that $\ker(L) \cap V= \{0\}$.
By a dimension argument, this imposes 
\be
n:= \dim(V) \le m.
\ee
In other words, we must place ourselves in an underparametrized regime for which the number $n$ of parameters describing the model does not exceed the number $m$ of data.
This contrasts with many current studies,
especially in the field of Deep Learning,
which emphasize the advantages of overparametrization.
In order to incorporate overparametrization in the optimal recovery problem under consideration,
we must then restrict the magnitude of model- and data-consistent elements.
A glaring strategy consists in altering the approximability set \eqref{ApproxMod}.
We do so in two different ways,
namely by considering a bounded approximability set {\em of the first type}, i.e.,
\be
\cK  = \cK^{\rm I}_{V,\eps,\kappa} :=
\big\{ f \in X: \; {\rm dist}_X(f,V) \le \eps \mbox{ and } \|f\|_X \le \kappa \big\}, 
\ee
and a bounded approximability set {\em of the second type}, i.e.,
\be
\cK = \cK^{\rm II}_{V,\eps,\kappa} := \big\{ f \in X: \; \exists v \in V \mbox{ with } \|f-v\|_X \le \eps \mbox{ and } \|v\|_X \le \kappa \big\}. 
\ee
We will start by analyzing the second type of bounded approximability sets in Section \ref{Sec2ndType} by formally describing the optimal recovery maps before revealing on a familiar example how the associated minimization problem is tackled in practice.
The main ingredient in essence belongs to the sum-of-squares techniques from semidefinite programming.
Next,
we will analyze the first type of bounded approximability sets in Section \ref{Sec1stType}.
We will even formally describe optimal recovery maps over more general model sets consisting of intersections of approximability sets.
On the prior example,
we will again reveal how the associated minimization problem  is tackled in practice.
This time, the main ingredient in essence belongs to the moment techniques from semidefinite programming.
In view of this article's emphasis on computability issues,
all of the theoretical constructions are illustrated in a reproducible {\sc matlab} file downloadable from the author's webpage.

\section{Bounded approximability set of the second type}\vspace{-2mm}
\label{Sec2ndType}

We concentrate in this section on the bounded approximability set of the second type,
i.e., on
\be
\label{TypeII}
\cK = \big\{ f \in X: \; \exists v \in V \mbox{ with } \|f-v\|_X \le \eps \mbox{ and } \|v\|_X \le \kappa \big\}.
\ee
We shall first describe optimal recovery maps before showing how they can be computed in practice.

\subsection{Description of an optimal recovery map}\vspace{-2mm}

The result below reveals how \cite[Theorem 3.1]{DFPW} extends from the model set~\eqref{ApproxMod}
to the model set~\eqref{TypeII}.

\bthm
\label{ThmDescrI}
If $Q: X \to \bR$ is a linear functional, then an optimal recovery map over the bounded approximability set \eqref{TypeII} is the linear functional
\be
R^{\rm opt}: y \in \bR^m \mapsto \sum_{i=1}^m a^{\rm opt}_i y_i \in \bR,
\ee
where the optimal weight $a^{\rm opt} \in \bR^m$ are precomputed as a solution to
\be
\label{OptTypeII}
\minimize{a \in \bR^m} 
\bigg[
\eps \times 
\Big\|
Q - \sum_{i=1}^m a_i \ell_i
\Big\|_{X^*} 
+  \kappa \times
\max_{v \in B_V} \Big|
Q(v) - \sum_{i=1}^m a_i \ell_i(v)
\Big| 
\bigg].
\ee
\ethm

\bpf
Since the model set \eqref{TypeII} is symmetric and convex,
there exists an optimal recovery map $R^{\rm opt}: \bR^m \to \bR$ which is linear,
i.e., of the form $R^{\rm opt}(y) = \sum_{i=1}^m a^{\rm opt}_i y_i$.
The vector $a^{\rm opt} \in \bR^m$ minimizes in particular the worst-case error
$e := \max \{ |Q(f) - \sum_{i=1}^m a_i \ell_i(f)| : f \in \cK \}$ among all $a \in \bR^m$.
Thus, it is sufficient to transform this worst-case error into the expression featured between square brackets in \eqref{OptTypeII}.
This is done by writing
\begin{align}
e & = \max_{f \in X}
\Big\{ \Big|Q(f) - \sum_{i=1}^m a_i \ell_i(f)\Big|:  \|f-v\|_X \le \eps \mbox{ for some } v \in V \mbox{ with } \|v\|_X \le \kappa \Big\}\\
\nonumber
& = \max_{\substack{f \in X\\ v \in V}}
\Big\{ \Big|Q(f) - \sum_{i=1}^m a_i \ell_i(f) \Big|: 
\|f-v\|_X \le \eps, \,
\|v\|_X \le \kappa \Big\}\\
\nonumber
& = \max_{\substack{h \in X\\ v \in V}}
\Big\{ \Big|Q(h+v) - \sum_{i=1}^m a_i \ell_i(h+v)\Big|: 
\|h\|_X \le \eps, \,
\|v\|_X \le \kappa \Big\}\\
\nonumber
& = \max_{h \in X}
\Big\{ \Big|Q(h) - \sum_{i=1}^m a_i \ell_i(h)\Big|: 
\|h\|_X \le \eps
\Big\}
+
\max_{v \in V}
\Big\{ \Big|Q(v) - \sum_{i=1}^m a_i \ell_i(v)\Big|: \|v\|_X \le \kappa \Big\}.
\end{align}
By homogeneity,
the latter is readily seen to coincide with the required expression.
\epf

\brk
The approximability set \eqref{ApproxMod} where the condition $\|v\|_X \le \kappa$ is not imposed can be viewed as an instantiation of \eqref{TypeII} with $\kappa = \infty$. 
In this instantiation,
if $\max_{v \in B_V} \Big|
Q(v) - \sum_{i=1}^m a_i \ell_i(v)
\Big|$ was nonzero,
then the objective function would be infinite.
Therefore, the infimum will be attained with the constraint $\max_{v \in B_V} \Big|
Q(v) - \sum_{i=1}^m a_i \ell_i(v)
\Big| = 0$ in effect.
This argument constitutes another way of deriving the form of the optimal recovery map over the original approximability set \eqref{ApproxMod}.
Let us note in passing that,
while the optimization program \eqref{ResultDFPW} 
was independent of $\eps>0$,
adding the condition $\|v\|_X \le \kappa$ does create  a dependence on $\eps > 0$
in the optimization program \eqref{OptTypeII}, unless $\kappa$ is proportional to~$\eps$.
\erk

\brk
In the presence of observation error $e \in \bR^m $ in $y = L(f)+ e$,
modeled as in \cite{EttFou} by the bounded uncertainty set
\be
\cE = \{ e \in \bR^m: \|e\|_p \le \eta \},
\ee
an optimal recovery map for a linear functional $Q: X \to \bR$
over $\cK$ and $\cE$ simultaneously still consists of a linear functional $R^{\rm opt}: y \in \bR^m \mapsto \sum_{i=1}^m a^{\rm opt}_i y_i$,
but now the optimal weights  $a^{\rm opt} \in \bR^m$ are solution to the optimization program
\be
\label{OptTypeIINoisy}
\minimize{a \in \bR^m} 
\bigg[
 \eps \times
\Big\|
Q - \sum_{i=1}^m a_i \ell_i
\Big\|_{X^*} 
+  \kappa \times
\max_{v \in B_V} \Big|
Q(v) - \sum_{i=1}^m a_i \ell_i(v)
\Big| 
 + \eta \times \|a\|_{p'}
\bigg],
\ee
where $p' = p/(p-1)$ is the conjugate exponent to $p \in [1,\infty]$.
The argument, which follows the ideas presented in \cite{EttFou},
is left to the reader.
We do point out that the program \eqref{OptTypeIINoisy} is solvable in practice as soon as soon as the program \eqref{OptTypeII} itself is solvable in practice,
for instance as in the forthcoming example.
\erk

\subsection{Computational realization for $X = \cC[-1,1]$}\vspace{-2mm}

For practical purposes, the result of Theorem \ref{ThmDescrI} is close to useless if the minimization cannot be performed efficiently.
We show below that in the important case $X = \cC[-1,1]$,
choosing $V$ as the space $\cP_n$ of algebraic polynomials of degree $<n$
leads to an optimization problem which can be solved exactly via semidefinite programming.
For that,
we also assume that the observation functionals are distinct point evaluations and that the quantity of interest $Q$ is another point evaluation or the normalized integral.
These restrictions can be lifted if we trade exact solutions for quantifiably approximate  solutions,
see Subsection \ref{SecCompI}.
In the statement below,
the notation ${\rm Toep}(x)$ represents the symmetric Toeplitz matrix built from a vector $x \in \bR^d$,
i.e.,
\be
{\rm Toep}(x) := \bbmx
x_1 & x_2 & \cdots & \cdots & x_d\\
x_2 & x_1 & x_2 & \ddots & \vdots \\
\vdots & \ddots& \ddots & \ddots & \vdots \\
\vdots & \ddots & x_2 & x_1 & x_2\\
 x_d & \cdots & \cdots & x_2 & x_1
\ebmx,
\ee
and the polynomials $T_j$, $j \in [0:n-1]$,
denote the $j$th Chebyshev polynomials of the first kind.
 
\bthm
Assuming that $V = \cP_n \inc \cC[-1,1]$ and that $\ell_1,\ldots,\ell_m$ are point evaluations at distinct points $x_1,\ldots,x_m \in [-1,1]$, 
an optimal recovery map over the bounded approximability set~\eqref{TypeII} for the quantity of interest $Q(f) = f(x_0)$, $x_0 \not\in \{ x_1,\ldots, x_m\}$, or $Q(f) = (1/2)\int_{-1}^1 f(x) dx$
is the linear functional
\be
R^{\rm opt}: y \in \bR^m \mapsto \sum_{i=1}^m a^{\rm opt}_i y_i \in \bR,
\ee
where the optimal weights $a^{\rm opt} \in \bR^m$ are precomputed as a solution to the semidefinite program
\begin{align}
\minimize{\substack{a,s \in \bR^m \\ u \in \bR^n}} \; \bigg[ \eps \times \sum_{i=1}^m s_i + \kappa \times u_1 \bigg]
& & \mbox{subject to } &  \; 
{\rm Toep}(u + C a - b) \succeq 0,
\;  {\rm Toep}(u - C a + b) \succeq 0,\\
\nonumber
& & \mbox{and } & \; 
s + a \ge 0, \; s-a \ge 0.
\end{align}
Here, $b \in \bR^n$ and $C \in \bR^{n \times m}$ have entries $b_j = Q(T_j)$ and $C_{j,i} = \ell_i(T_j)$, $i \in [1:m]$, $j \in [0:n-1]$.
\ethm

\bpf
The work consists in recasting the objective function of \eqref{OptTypeII}
into manageable form.
Under the assumptions on $\ell_1,\ldots,\ell_m$ and on $Q$,
the first term is not a problem, by virtue of
\be
\Big\| Q  - \sum_{i=1}^m a_i \ell_i \Big\|_{\cC[-1,1]^*}  = 1 + \sum_{i=1}^m |a_i|.
\ee
We now turn to the second term,
i.e., the one involving the  maximum over the unit ball $B_V$ of~$V$.
The idea, common in Robust Optimization \cite{RO}, 
relies on duality to change the maximum into a minimum,
which is then integrated into a larger  minimization problem. 
This is possible essentially when $B_V$ admits a linear or semidefinite description,
which is the case for $V = \cP_n$.
Indeed, as already observed in \cite[Subsection 5.3]{MinProjExt},
following ideas formulated in \cite{FouPow},
the unit ball of $\cP_n$ admits the semidefinite description
\begin{align}
B_{\cP_n} = \bigg\{
\sum_{j=0}^{n-1} \tr [D_j(P-M)] T_j
& \mbox{ for some positive semidefinite matrices }M,P \in \bR^{n \times n} \\
\nonumber
& \mbox{ that satisfy }
 \tr[D_j(P+M)]
= \delta_{0,j}
\bigg\},
\end{align}
where, for each $j \in [0:n-1]$,
the symmetric matrix
\be
D_j := \bbmx
0 & \cdots & 0 & 1  & 0 & \cdots & 0\\
\vdots & 0 & \ddots & 0 & 1 & \ddots & \vdots \\
0 & \ddots & \ddots & \ddots & \ddots & \ddots & 0\\
1 & 0 & \ddots & \ddots & \ddots & 0 & 1\\
0 & 1 & \ddots & \ddots & \ddots & \ddots & 0\\
\vdots & \ddots & \ddots & 0 & \ddots & 0 & \vdots\\
0 & \cdots & 0 & 1 & 0 & \cdots & 0
\ebmx \in \bR^{n \times n}
\ee
has $1$'s on the $j$th subdiagonal and superdiagonal and $0$'s elsewhere --- in particular $D_0$ is the $n \times n$ identity matrix.
Thus, for a fixed $a \in \bR^m$,
with $Q_a := Q - \sum_{i=1}^m a_i \ell_i$,
the maximum over $B_V$ reads
\be
\max_{M,P \in \bR^{n \times n}} \bigg\{ \tr \Big[ \Big(
\sum_{j=0}^{n-1} Q_a(T_j) D_j \Big)
 (P-M)\Big]
 : M, P \succeq 0, \,
\tr[D_j(P+M)]
= \delta_{0,j}
\bigg\}.
\ee
Invoking duality in semidefinite programming
(see e.g. \cite[p.265-266]{BoyVan}), the latter can be transformed into
\be
\label{above}
\min_{u \in \bR^n} u_1 
\qquad \mbox{subject to }
\sum_{j=0}^{n-1} u_j D_j
\pm \sum_{j=0}^{n-1} Q_a(T_j) D_j
\succeq 0.
\ee
Since $Q_a(T_j) = b_j - (Ca)_j$ for any $j \in [0:n-1]$, the constraint in \eqref{above} can be condensed to ${\rm Toep}( u \pm (Ca - b) ) \succeq 0$.
Then, combining the minimization over $u \in \bR^n$ with the minimization over $a \in \bR^m$,
the optimization program \eqref{OptTypeII} becomes equivalent to 
\be
\minimize{\substack{a \in \bR^m\\ u \in \bR^n}} \; \bigg[ \eps \times \sum_{i=1}^m |a_i| + \kappa \times u_1 \bigg]
\qquad \mbox{subject to }
{\rm Toep}( u \pm (Ca - b) ) \succeq 0.
\ee
The final step is the introduction of  slack variables $s \in \bR^m$ such that $|a_i| \le s_i$,
i.e., $-s_i \le a_i \le s_i$, 
for all $i \in [1:m]$.
\epf

\section{Bounded approximability set of the first type}\vspace{-2mm}
\label{Sec1stType}

We concentrate in this section on the bounded approximability set of the first type, i.e., on
\be
\label{TypeI}
\cK = \big\{ f \in X: \; {\rm dist}_X(f,V) \le \eps \mbox{ and } \|f\|_X \le \kappa \big\}. 
\ee
Once again,
we shall first describe optimal recovery maps before showing how they can be computed in practice.

\subsection{Description of an optimal recovery map}\vspace{-2mm}

The result below reveals how \cite[Theorem 3.1]{DFPW} extends from the model set \eqref{ApproxMod} to the model set~\eqref{TypeI}.

\bthm
\label{ThmTypeI}
If $Q: X \to \bR$ is a linear functional, then an optimal recovery map over the bounded approximability set \eqref{TypeI} is the linear functional
\be
\label{RoptTypeI}
R^{\rm opt}: y \in \bR^m \mapsto \sum_{i=1}^m a^{\rm opt}_i y_i \in \bR,
\ee
where the optimal weights $a^{\rm opt} \in \bR^m$ are precomputed as a solution to
\be
\label{OptTypeI}
\minimize{\substack{a \in \bR^m\\ \mu,\nu \in X^* }} 
\; \big[ \eps \times \|\mu\|_{X^*} + \kappa \times \|\nu\|_{X^*} \big]
\qquad \mbox{subject to } \mu + \nu = Q-\sum_{i=1}^m a_i \ell_i
\; \mbox{ and } \; \mu_{|V} = 0.
\ee
\ethm

As a matter of fact,
Theorem \ref{ThmTypeI} is a corollary of Theorem \ref{ThmGen} below. 
The setting of the more general result involves subspaces $V_1,\ldots,V_K$ of a linear space $X$ equipped with possibly distinct norms $\|\cdot\|_{(1)},\ldots,\|\cdot\|_{(K)}$.
The model set is then defined, 
for some parameters $\eps_1,\ldots,\eps_K >0$,
by
\be
\label{Inters}
\cK = \{ f \in X: {\rm dist}_{\|\cdot\|_{(1)}}(f,V_1) \le \eps_1,
\ldots,  {\rm dist}_{\|\cdot\|_{(K)}}(f,V_K) \le \eps_K \}.
\ee
It corresponds to what was called the multispace problem in \cite[Section 3]{BCDDPW}.
One works under the assumption that
\be
\label{Assump}
\ker(L) \cap V_1 \cap \ldots \cap V_K  = \{ 0 \}.
\ee
This assumption holds for the bounded approximability set of the first type,
obtained by taking $V_1 = V$, $V_2 = \{0\}$, and $\|\cdot\|_{(1)} = \|\cdot\|_{(2)} = \|\cdot\|_{X}$.

\bthm
\label{ThmGen}
If $Q: X \to \bR$ is a linear functional, then an optimal recovery map over the model set \eqref{Inters} is the linear functional
\be
R^{\rm opt}:
y \in \bR^m \mapsto \sum_{i=1}^m a^{\rm opt}_i y_i \in \bR,
\ee
where the optimal weights $a^{\rm opt} \in \bR^m $ are precomputed as a solution to
\begin{align}
\label{OptTypeIGen}
\minimize{\substack{a \in \bR^m\\
\la_1,\ldots,\la_K \in X^*}} \;
\big[ \eps_1 \|\la_1\|_{(1)}^* + \cdots  + \eps_K \|\la_K\|_{(K)}^* \big]
& & \mbox{subject to } &
 \la_1 + \cdots + \la_K = Q - \sum_{i=1}^m a_i \ell_i\\
 \nonumber
& & \mbox{and } &
 {\la_1}_{|V_1} = 0, \ldots, {\la_K}_{|V_K} = 0.
\end{align}
\ethm

\bpf
We first notice that, replacing the norms $\|\cdot\|_{(k)}$ by $\|\cdot\|_{(k)}/\eps_k$,
we can assume that $\eps_k = 1$. 
Next, since the model set $\cK$ is symmetric and convex,
there exists an optimal recovery map which is linear,
i.e., of the form $ y \in \bR^m \mapsto \sum_{i=1}^m a_i y_i \in \bR$.
An optimal weight vector $a^{\rm opt} \in \bR^m$ is then obtained as a solution to the optimization problem
\be 
\label{APriori}
\minimize{a \in \bR^m}
\max_{f \in X} \bigg\{
\bigg| Q(f) - \sum_{i=1}^m a_i \ell_i(f) \bigg|: \; {\rm dist}_{\|\cdot\|_{(k)}}(f,V_k) \le 1  \mbox{ for all } k \in [1:K]
\bigg\}.
\ee
We claim that an optimal weight vector $a^{\rm opt} \in \bR^m$ can also be obtained as a solution to the optimization problem
\begin{align}
\label{Implem}
\minimize{a \in \bR^m}
\min_{ \la_1,\ldots,\la_K \in X^*}
\bigg\{
\|\la_1\|_{(1)}^* + \cdots  + \|\la_K\|_{(K)}^*: & \;
\la_1 + \cdots + \la_K = Q - \sum_{i=1}^m a_i \ell_i,\\
\nonumber
& \;  {\la_k}_{|V_k} =0 \mbox{ for all } k \in [1:K]
\bigg\}.
\end{align}
In other words, we shall prove in two steps that the minimal values of \eqref{APriori} and \eqref{Implem} coincide.

Firstly, we shall justify that the objective function in \eqref{APriori} is bounded by the objective function in~\eqref{Implem}
--- a property which holds independently of $a \in \bR^m$.
To do so, let us consider $f \in X$ such that $\|f-v_1\|_{(1)}\le 1, \ldots, \|f-v_K\|_{(K)} \le 1$ for some $v_1 \in V_1$, $\ldots$, $v_K \in V_K$.
Let us also consider $\la_1,\ldots,\la_K \in X^*$ such that $\la_1 + \cdots + \la_K = Q - \sum_{i=1}^m a_i \ell_i$ and ${\la_1}_{|V_1} =0$, $\ldots$, ${\la_K}_{|V_K} =0$.
We have
\begin{align}
\bigg| Q(f) - \sum_{i=1}^m a_i \ell_i(f) \bigg|
& = |\la_1(f) + \cdots + \la_K(f)|
= |\la_1(f-v_1) + \cdots + \la_K(f-v_K)|\\
\nonumber
& \le \|\la_1\|^*_{(1)} \|f-v_1\|_{(1)} + \cdots  +  \|\la_K\|^*_{(1)} \|f-v_K\|_{(K)}\\
\nonumber
& \le \|\la_1\|_{(1)}^* + \cdots  + \|\la_K\|_{(K)}^*.
\end{align}
Taking the infimum over $\la_1,\ldots,\la_K$ and the supremum over $f$ yields the desired result.

Secondly,
we shall justify that the minimal value of~\eqref{Implem} is bounded by the minimal value of \eqref{APriori}.
To do so, let us consider the linear space $Z:= X \times \cdots \times X$
equipped with the norm
\be
\|(f_1,\ldots,f_K)\|_Z := \max_{k \in [1:K]} \|f_k\|_{(k)}. 
\ee
Introducing the subspace $U$ of $Z$ given by 
\be
U := \{ (h,\ldots,h), h \in \ker(L)\}, 
\ee
the assumption \eqref{Assump} is equivalent to $U \cap (V_1 \times \cdots \times V_K) = \{ 0 \}$.
Thus, we can define a linear functional $\la$ on $U \oplus (V_1 \times \cdots \times V_K) $ by
\begin{align}
\la((h,\ldots,h )) & = Q(h)
& & \mbox{ for } h \in \ker(L),\\
\la( (v_1,\ldots,v_K) ) & = 0
& & \mbox{ for } (v_1,\ldots,v_K) \in V_1 \times \cdots \times V_K.
\end{align}
Let then $\wt{\la} \in Z^*$ denote a Hahn--Banach extension of $\la$ to the whole $Z$.
With linear functionals $\la_1,\ldots,\la_K \in X^*$
defined for each $k \in [1:K]$ and $f \in X$ by $\la_k(f) = \wt{\la}((0,\ldots,0,f,0, \ldots,0))$, where $f$ appears at the $k$th position,
we have $\la_1(f) + \cdots + \la_K(f) = \wt{\la}((f,\ldots,f))$ for all $f \in X$,
hence in particular $Q - (\la_1+\cdots+\la_K)$ vanishes on $\ker(L)$.
This implies (see e.g. \cite[Lemma 3.9]{Rud}) that $Q - (\la_1+\cdots+\la_K) = \sum_{i=1}^m a^\sharp_i \ell_i$ for some $a^\sharp \in \bR^m$.
In other words, the first constraint
in~\eqref{Implem}
is satisfied by $a^\sharp$ and $\la_1,\ldots,\la_K$.
The second constraint is also satisfied: 
indeed, for $v_k \in V_k$,  $\la_k(v_k) = \wt{\la}((0,\ldots,0,v_k,0,\ldots,0)) =0$ since $(0,\ldots,0,v_k,0,\ldots,0) \in V_1 \times \cdots \times V_K$.
Therefore, the minimal value of \eqref{Implem} is bounded by
\begin{align}
\|\la_1\|_{(1)}^* + \cdots  + \|\la_K\|_{(K)}^*
& = 
\max_{\|f_1\|_{(1)} \le 1} \la_1(f_1) 
+ \cdots + 
\max_{\|f_K\|_{(K)}\le 1}  \la_K(f_K) \\
\nonumber
& = \max_{\|f_1\|_{(1)} \le 1, \ldots, \|f_K\|_{(K)}\le 1} \la_1(f_1) + \cdots + \la_K(f_K)\\
\nonumber
& = \max_{\|(f_1,\ldots,f_K)\|_Z \le 1} \wt{\la}((f_1,\ldots,f_K))\\
\nonumber
& = \|\wt{\la}\|_{Z}^*.
\end{align}
The latter equals the norm of $\la$ on $U \oplus (V_1 \times \cdots \times V_K)$,
by virtue of $\wt{\la}$ being a Hahn--Banach extension of $\la$, so that
\begin{align}
\|\la_1\|_{(1)}^* + \cdots  + \|\la_K\|_{(K)}^*
& = \max_{\substack{u = (h,\ldots,h) \in U\\
v = (v_1,\ldots,v_K) \in V_1 \times \cdots \times V_K}}
\big\{
\la(u-v):
\|u-v\|_Z \le 1
\big\}\\
\nonumber
& = \max_{\substack{h \in \ker(L) \\ v_k \in V_k}}
\big\{
Q(h) : \|h-v_k\|_{(k)} \le 1 \mbox{ for all } k \in [1:K]
\big\}.
\end{align}
It follows that, for any $a \in \bR^m$, 
\begin{align}
\|\la_1\|_{(1)}^* + \cdots  + \|\la_K\|_{(K)}^* 
& = \max_{\substack{h \in \ker(L)\\
v_k \in V_k }}
\bigg\{
Q(h) - \sum_{i=1}^m a_i \ell_i(h):  \|h-v_k\|_{(k)} \le 1 \mbox{ for all } k \in [1:K]
\bigg\}\\
\nonumber
& \le
\max_{f \in X} \bigg\{
\bigg| Q(f) \hspace{-1mm} - \hspace{-1mm} \sum_{i=1}^m a_i \ell_i(f) \bigg|: {\rm dist}_{\|\cdot\|_{(k)}}(f,V_k) \le 1 \mbox{ for all }k \in [1\hspace{-1mm}:\hspace{-1mm}K]
\bigg\}.
\end{align}
Taking the minimum over all $a \in \bR^m$ shows that$\|\la_1\|_{(1)}^* + \cdots  + \|\la_K\|_{(K)}^*$
is less than or equal to the minimal value of \eqref{APriori},
and in turn that the same is true for
the minimal value of \eqref{Implem}.
\epf 

\brk
The approximability set \eqref{ApproxMod} where the condition $\|f\|_X  \le \kappa$ is not imposed can be viewed as an instantiation of \eqref{TypeI} with  $\kappa = \infty$. 
In this instantiation,
if $\| \nu \|_{X^*}$ was nonzero,
then the objective function in \eqref{OptTypeI} would be infinite.
Therefore, the minimum will be attained with the constraint $\| \nu \|_{X^*} = 0$ in effect,
leading to $\mu = Q - \sum_{i=1}^m a_i \ell_i$
and in turn to the constraint $(Q - \sum_{i=1}^m a_i \ell_i)_{|V} = 0$.
We do retrieve the minimization of  \eqref{ResultDFPW}, as expected.
We note in passing that,
while the optimization program \eqref{ResultDFPW} was independent of $\eps>0$,
adding the condition $\|f\|_X \le \kappa$ does create a dependence on $\eps > 0$ in the optimization problem \eqref{OptTypeI}, unless $\kappa$ is proportional to $\eps$.
\erk

\brk
In the presence of observation error $e \in \bR^m $ in $y = L(f)+ e$,
again modeled as in \cite{EttFou} by the bounded uncertainty set
\be
\cE = \{ e \in \bR^m: \|e\|_p \le \eta \},
\ee
an optimal recovery map for a linear functional $Q: X \to \bR$
over $\cK$ and $\cE$ simultaneously still consists of a linear functional $R^{\rm opt}: y \in \bR^m \mapsto \sum_{i=1}^m a^{\rm opt}_i y_i$,
but now the optimal weights $a^{\rm opt} \in \bR^m$ are solution to the optimization program
\be
\label{OptTypeINoisy}
\minimize{\substack{a \in \bR^m\\ \mu,\nu \in X^* }} 
\; \big[ \eps \times \|\mu\|_{X^*} + \kappa \times \|\nu\|_{X^*} +
\eta \times \|a\|_{p'} \big]
\qquad \mbox{subject to } \mu + \nu = Q-\sum_{i=1}^m a_i \ell_i
\; \mbox{ and } \; \mu_{|V} = 0.
\ee
The argument follows the ideas presented in \cite{EttFou}
and, although more subtle,
is once again left to the reader.
We do point out that the program \eqref{OptTypeINoisy} is solvable in practice as soon as soon as the program \eqref{OptTypeI} itself is solvable in practice, for instance as in the forthcoming example.
\erk

\subsection{Computational realization for $X = \cC[-1,1]$}\vspace{-2mm}
\label{SecCompI}

As before,
the high-level results of Theorems \ref{ThmTypeI} and \ref{ThmGen} are of little practical use if the minimizations \eqref{OptTypeI} and \eqref{OptTypeIGen} cannot be performed efficiently.
In the important situation $X = \cC[-1,1]$,
the dual functionals appearing as optimization variables are identified with measures.
Despite involving infinite dimensional objects,
minimizations over measures can be tackled via semidefinite programming, see e.g. \cite{Las}.
Although such minimizations are in general not solved exactly,
their accuracy can be quantifiably estimated in our specific case.
For ease of presentation, we illustrate the approach by concentrating on the optimization program \eqref{OptTypeI} rather than \eqref{OptTypeIGen}.
We also assume that $V = \cP_n$ 
and we write the observation functionals $\ell_1,\ldots,\ell_m$, 
as well as the quantity of interest $Q$,
as
\be
\ell_i(f) = \int_{-1}^1 f(x) d\la_i(x),
\qquad
Q(f) = \int_{-1}^1 f(x) d\rho(x),
\qquad f \in \cC[-1,1],
\ee
for some signed Borel measures $\la_1,\ldots,\la_m,\rho$ defined on $[-1,1]$.
In this way, passing from linear functionals to signed Borel measures as optimization variables,
the program \eqref{OptTypeI} reads
\begin{align}
\label{OptTypeIonC}
\minimize{\substack{a \in \bR^m\\ \mu,\nu }} 
\;  \int_{-1}^1 \eps \, d|\mu| + \, \kappa \, d|\nu|
& &  \mbox{subject to }& \mu + \nu = \rho - \sum_{i=1}^m a_i \la_i \\
\nonumber
& & \mbox{ and }& \int_{-1}^1 v(x) d\mu(x) = 0 \mbox{ for all } v \in \cP_n.
\end{align}
Let us introduce as slack variables the nonnegative Borel measures $\mu^+$, $\mu^-$, $\nu^+$, and $\nu^-$ involved in the Jordan decompositions $\mu = \mu^+ - \mu^-$ and $\nu = \nu^+ - \nu^-$,
so that the problem \eqref{OptTypeIonC}
is recast as
\begin{align}
\label{OptTypeIonC2}
\minimize{\substack{a \in \bR^m\\ \mu^\pm, \nu^\pm}}
\int_{-1}^1 \eps \, d(\mu^+ + \mu^-) + \kappa \, d(\nu^+ + \nu^-)
& & \mbox{s.to } & 
\mu^+-\mu^- + \nu^+ - \nu^- = \rho -\sum_{i=1}^m a_i \la_i\\
\nonumber
& & \mbox{and } & 
\int_{-1}^1 v(x) d(\mu^+-\mu^-)(x) = 0
\mbox{ for all } v \in \cP_n.
\end{align}
Next, replacing the measures $\mu^\pm$ and $\nu^\pm$ by the infinite sequences of moments $w^\pm = \cM_\infty(\mu^\pm) \in \bR^\bN$ and $z^\pm = \cM_\infty(\nu^\pm) \in \bR^\bN$ of moments defined for $k \ge 1$ by
\be
w^\pm_k = 
\int_{-1}^1 T_{k-1}(x) d\mu^\pm(x),
\qquad \quad 
z^\pm_k = 
\int_{-1}^1  T_{k-1}(x) d \nu^\pm(x),
\ee
the problem \eqref{OptTypeIonC2} is equivalent\footnote{the equivalence is based on the discrete trigonometric moment problem, see \cite{FouLas} for details.} to the infinite semidefinite program
\begin{align}
\label{OptTypeIonC3}
\minimize{\substack{a \in \bR^m \\
w^\pm,z^\pm \in \bR^\bN}}
\; \eps \, (w_1^+ + w_1^-) + \kappa \, (z_1^+ + z_1^-),
& & \mbox{s.to } & w^+ - w^- + z^+ - z^- = \cM_\infty \Big(\rho-\sum_{i=1}^m a_i \la_i \Big),\\
\nonumber
& & \mbox{and } & w^+_j - w^-_j = 0 \mbox{ for all } j \in [1:n],\\ 
\nonumber
& & \mbox{and } & {\rm Toep}_\infty(w^\pm) \succeq 0, \; 
{\rm Toep}_\infty(z^\pm) \succeq 0.
 \end{align}
Instead of solving this infinite optimization program,
we truncate it to a level $N \ge n$ and solve instead the resulting finite semidefinite program
\begin{align}
\label{OptTypeITrunc}
\minimize{\substack{a \in \bR^m \\
w^\pm,z^\pm \in \bR^N}}
\; \eps \, (w_1^+ + w_1^-) + \kappa \, (z_1^+ + z_1^-),
& & \mbox{s.to } & w^+ - w^- + z^+ - z^- = \cM_N \Big(\rho-\sum_{i=1}^m a_i \la_i \Big)\\
\nonumber
& & \mbox{and } & w^+_j - w^-_j = 0 \mbox{ for all } j \in [1:n],\\ 
\nonumber
& & \mbox{and } & {\rm Toep}_N(w^\pm) \succeq 0, \; 
{\rm Toep}_N(z^\pm) \succeq 0.
 \end{align}
The rest of this section is devoted to justifying in a  quantitative way that the minimal value of this truncated problem converges to the minimal value of the original problem.
We also justify,
although not quantitatively, 
 that the vectors $a^{(N)} \in \bR^m$ obtained by solving \eqref{OptTypeITrunc} converge as $N \to \infty$ to a vector $a^{\rm opt} \in \bR^m$  solving \eqref{OptTypeI}
--- we do not analyze the behavior of the optimal measures since only the vector $a^{\rm opt}$ is required in \eqref{RoptTypeI}.
From now on,
we also work under the assumption of linear independence for the restrictions ${\ell_1}_{|\cP_m}, \ldots, {\ell_m}_{|\cP_m}$ of the observation functionals to the space $\cP_m$ of polynomials of degree $<m$.
This assumption is easily seen to be equivalent to the invertibility of (the transpose of) the moment matrix $M \in \bR^{m \times m}$ defined by
\be
M_{j,i} = \ell_i(T_{j-1}),
\qquad i,j \in [1:m].
\ee
This holds e.g. when the observation functionals are evaluations at $m$ distinct points in $[-1,1]$.

\bthm
\label{ThmaN2aopt}
Suppose that the system $( {\ell_1}_{|\cP_m}, \ldots, {\ell_m}_{|\cP_m} )$ is linearly independent.
If there is a unique $a^{\rm opt} \in \bR^m$ yielding  a minimizer $(a^{\rm opt},\mu^{\rm opt}, \nu^{\rm opt})$ of \eqref{OptTypeI},
then $a^{\rm opt}$ is the limit of any sequence $(a^{(N)})_{N \ge n}$ obtained by solving \eqref{OptTypeITrunc} for each $N \ge n$.
Without uniqueness,
it still holds that any subsequence of $(a^{(N)})_{N\ge n}$ admits a subsequence converging to the first component of a minimizer of \eqref{OptTypeI}.
\ethm

\bpf
The first part of the theorem follows from the second part:
it is indeed well-known that the convergence of a sequence towards a given point is guaranteed as soon as any of its subsequences admits a subsequence converging to that point.

To establish the second part,
let $\alpha^{(N)} \in \bR$ and  $(a^{(N)},w^{\pm,(N)},z^{\pm,(N)}) \in \bR^m \times (\bR^{N})^4$ denote, for each $N \ge n$, the minimum value and some minimizer of \eqref{OptTypeITrunc}, respectively.
We write $w^{\pm,((N))} \in \bR^\bN$ and $z^{\pm,((N))} \in \bR^\bN$ for the infinite vectors obtained by padding the finite vectors $w^{\pm,(N)} \in \bR^N$ and  $z^{\pm,(N)} \in \bR^N$ with zeros.
Let us now consider a subsequence $\big( (a^{(N_k)},w^{\pm,((N_k))},z^{\pm,((N_k))}) \big)_{k \ge 1}$ of the $(\ell_\infty^m \times (\ell_\infty^\bN)^4)$-valued sequence $\big( (a^{(N)},w^{\pm,(N)},z^{\pm,(N)})\big)_{N \ge n}$.
Our objective is to show that there exist a subsequence $\big( (a^{(N_{k_\ell})},w^{\pm,((N_{k_\ell}))},z^{\pm,((N_{k_\ell}))}) \big)_{\ell \ge 1}$
and a minimizer $(\wt{a},\wt{w}^{\pm},\wt{z}^{\pm})$ of \eqref{OptTypeIonC3}
such that $a^{(N_{k_\ell})}$ converges to $\wt{a}$ as $\ell \to \infty$.
To this end, we start by observing that the sequence $(\alpha^{(N_k)})_{k \ge 1}$ is nondecreasing and bounded by the minimal value $\alpha^{\rm opt}$ of \eqref{OptTypeIonC3}:
firstly, 
the inequality $\alpha^{(N_k)} \le \alpha^{(N_{k+1})}$ 
follows from the feasibility of $(a^{(N_{k+1})},w^{\pm,(N_{k+1})}_{[1:N_k]},z^{\pm,(N_{k+1})}_{[1:N_k]})$
for \eqref{OptTypeITrunc} specified to $N=N_k$, so that 
\be
\alpha^{(N_k)}
\le \eps \, (w_1^{+,(N_{k+1})} + w_1^{-,(N_{k+1})}) + \kappa \, (z_1^{+,(N_{k+1})} + z_1^{-,(N_{k+1})}) = \alpha^{(N_{k+1})};
\ee
secondly, the inequality $\alpha^{(N_k)} \le \alpha^{\rm opt}$ similarly follows from the feasibility of $(a^{\rm opt},w^{\pm ,\rm opt}_{[1:N_k]},z^{\pm ,\rm opt}_{[1:N_k]})$ for \eqref{OptTypeITrunc} specified to $N=N_k$,
where evidently $(a^{\rm opt},w^{\pm ,\rm opt},z^{\pm ,\rm opt})$ represents some minimizer of \eqref{OptTypeIonC3}.
We continue by remarking that the positive semidefiniteness of ${\rm Toep}_{N_k}(w^{\pm,(N_k)})$ and ${\rm Toep}_{N_k}(z^{\pm,(N_k)})$ yields,
 for any $j \in [1:N_k]$,
\begin{align}
\big| w^{\pm,(N_k)}_j \big|  
& \le w^{\pm,(N_k)}_1
\le \f{1}{\eps} \big( 
\eps \, (w_1^{+,(N_{k})} + w_1^{-,(N_{k})}) + \kappa \, (z_1^{+,(N_{k})} + z_1^{-,(N_{k})})
\big) 
\le \f{1}{\eps} \alpha^{(N_k)}
\le \f{1}{\eps} \alpha^{\rm opt},
\\
\big| z^{\pm,(N_k)}_j \big|  
& \le z^{\pm,(N_k)}_1
\le \f{1}{\kappa} \big( 
\eps \, (w_1^{+,(N_{k})} + w_1^{-,(N_{k})}) + \kappa \, (z_1^{+,(N_{k})} + z_1^{-,(N_{k})})
\big) 
\le \f{1}{\kappa} \alpha^{(N_k)}
\le \f{1}{\kappa} \alpha^{\rm opt}.
\end{align}
Thus,
the $\ell_\infty^\bN$-valued sequences
$(w^{\pm,((N_k))})_{k \ge 1}$ and
$(z^{\pm,((N_k))})_{k \ge 1}$
are bounded.
This guarantees, 
by the sequential compactness Banach--Alaoglu theorem,
that $(w^{\pm,((N_k))})_{k \ge 1}$ and  $(z^{\pm,((N_k))})_{k \ge 1}$ admit convergent subsequences in the weak-star topology of $\ell_\infty^\bN$.
We denote the resulting convergent subsequence and its limit by
$\big( (w^{\pm,((N_{k_\ell}))},z^{\pm,((N_{k_\ell}))}) \big)_{\ell \ge 1}$ 
and $(\wt{w}^\pm,\wt{z}^\pm)$, respectively.
They come with associated $(a^{(N_{k_\ell})})_{\ell \ge 1}$ and $\wt{a}$.
Indeed,
the constraint imposed on minimizers of \eqref{OptTypeITrunc} yields in particular,
for $N_{k_\ell} \ge m$,
\be
M a^{(N_{k_\ell})} = \cM_m(\rho) 
- ( w^{+, (N_{k_\ell})} - w^{-, (N_{k_\ell})} + z^{+, (N_{k_\ell})} - z^{-, (N_{k_\ell})})_{[1:m]}.
\ee
In view of $w^{\pm,(N_{k_\ell})}_j \to \wt{w}^{\pm}_j$ and
$z^{\pm,(N_{k_\ell})}_j \to \wt{z}^{\pm}_j$ for all $j \ge 1$, 
which is a consequence of the weak-star convergence,
and of the invertibility assumption for the moment matrix,
we see that
\be
a^{(N_{k_\ell})} \underset{\ell \to \infty}{\tto} \wt{a} := M^{-1}  
\big(
\cM_m(\rho) 
- ( \wt{w}^{+} - \wt{w}^{-} + \wt{z}^{+} - \wt{z}^{-})_{[1:m]}
\big).
\ee
It remains to prove that the quintuple $(\wt{a},\wt{w}^\pm, \wt{z}^\pm)$ is a minimizer of \eqref{OptTypeIonC3}.
Writing the constraint of \eqref{OptTypeITrunc} satisfied by $(a^{(N_{k_\ell})},w^{\pm,(N_{k_\ell})},z^{\pm,(N_{k_\ell})})$
 and passing to the limit as $\ell \to \infty$
shows that the quintuple is feasible for \eqref{OptTypeIonC3}.
It is also a minimizer for this program,
by virtue of
\begin{align}
\eps \, (\wt{w}_1^{+} + \wt{w}_1^{-}) + \kappa \, (\wt{z}_1^{+} + \wt{z}_1^{-})
& = \lim_{\ell \to \infty}
\big(
\eps \, (w_1^{+,(N_{k_\ell})} + w_1^{-,(N_{k_\ell})}) + \kappa \, (z_1^{+,(N_{k_\ell})} + z_1^{-,(N_{k_\ell})})
\big) \\
\nonumber
& = \lim_{\ell \to \infty} \alpha^{(N_{k_\ell})} 
\le \alpha^{\rm opt}.
\end{align}
Our objective is now established, 
so the second part of theorem is proved.
\epf

Theorem \ref{ThmaN2aopt} does not tell us how to choose $N$ in order to reach a prescribed accuracy on $\|a^{\rm opt} - a^{(N)}\|$, not even on $\alpha^{\rm opt} - \alpha^{(N)}$.
We intend to provide an accuracy estimate for the latter,
which we do under the further restriction that the observation functionals and the quantity of interest are point evaluations at distinct $x_1,\ldots,x_m \in [-1,1]$
and at $x_0 \in [-1,1] \setm \{x_1,\ldots,x_m\}$, respectively.
Moreover, the estimate is an {\sl a posteriori} one,
in the sense that we need to solve problem \eqref{OptTypeITrunc} for a particular $N$ first.
As a matter of fact, we also need to solve an extra linear program subordinated to a particular grid $\ul{t} = (t_1,\ldots,t_K)$ avoiding the $x_i$'s.
This program is 
\begin{align}
\label{LPforUB}
\minimize{\substack{a \in \bR^m\\ u,v \in \bR^{1+m+K}\\ r,s \in \bR^{1+m+K}}}
\eps \sum_{h} r_h
+ \kappa \sum_{h} s_h
& &  \mbox{s.to } & 
u+v = \bbmx 1\\ \hline -a \\
\hline 0\ebmx,
\, \bbmx \, b \, | \, C \, | \, D \, \ebmx u = 0,\\
\nonumber
& & \mbox{and } & 
r + u \ge 0, \, r - u \ge 0,
s + v \ge 0, \, s - v \ge 0.
\end{align}
Here, the vector $b \in \bR^n$, the matrix $C \in \bR^{n \times m}$ (both encountered before),
and the matrix $D \in \bR^{n \times K}$ have entries
$b_j = T_{j-1}(x_0)$,
$C_{j,i} = T_{j-1}(x_i)$, and $D_{j,k} = T_{j-1}(t_k)$, $i \in [1:m]$, $j \in [1:n]$, $k \in [1:K]$.

\bthm
Suppose that the observation functionals $\ell_1,\ldots,\ell_m$ and the quantity of interest $Q$ take the form $\ell_i(f) = f(x_i)$, $i \in [1:m]$,
and $Q(f) = f(x_0)$
for distinct points $x_0, x_1,\ldots,x_m \in [-1,1]$.
For any $N \ge n$ and any grid $\ul{t}$ of $[-1,1]$ not intersecting $\{ x_0,x_1,\ldots,x_m\}$,
the minimal value $\alpha^{\rm opt}$ of \eqref{OptTypeI} satisfies
\be
\label{APostEst}
\alpha^{(N)} \le \alpha^{\rm opt} \le \beta^{(\ul{t})},
\ee
where $\alpha^{(N)}$ is the minimal value of \eqref{OptTypeITrunc} and $\beta^{(\ul{t})}$ is the minimal value of \eqref{LPforUB}. 
\ethm

\bpf
The leftmost inequality of \eqref{APostEst} was already justified (implicitly) in the proof of Theorem~\ref{ThmaN2aopt}
and actually does not rely on any assumption on the observation functionals or the quantity of interest.
For the rightmost inequality of \eqref{APostEst},
we keep in mind that $\alpha^{\rm opt}$ is the minimal value of the optimization program \eqref{OptTypeIonC}.
Since this is a minimal value over all signed Borel measures,
an upper bound is provided by the minimal value $\beta^{(\ul{t})}$ over the subset of all signed Borel measures consisting of linear combinations of Dirac measures at the distinct points $x_0,x_1,\ldots,x_m,t_1,\ldots,t_K$
--- the points $x_0,x_1,\ldots,x_m$ are included in order to make the constraint of \eqref{OptTypeIonC} feasible.
Writing such measures as
\begin{align}
\mu & = u \, \delta_{x_0} + \sum_{i=1}^m u'_i \delta_{x_i} + \sum_{k=1}^K u''_k \delta_{t_k},\\
\nu & = v \, \delta_{x_0} + \sum_{i=1}^m v'_i \delta_{x_i} + \sum_{k=1}^K v''_k \delta_{t_k},
\end{align}
we see that the upper bound $\beta^{(\ul{t})}$ takes the form
\begin{align}
\beta^{(\ul{t})}
 = \underset{\substack{a \in \bR^m, u,v \in \bR\\ u',v' \in \bR^m, u'',v'' \in \bR^K}}{\min}
\Big\{ 
\eps  (|u| +& \|u'\|_1 + \|u''\|_1) 
+ 
\kappa  (|v| + \|v'\|_1 + \|v''\|_1) :
\\
\nonumber
&  u+v = 1, \, u'+v' = -a, \, u''+v''=0, \, u b + C u' + D u'' =0
\Big\}.
\end{align}
By gathering $u \in \bR$, $u' \in \bR^m$,
and $u'' \in \bR^K$ into a single vector $\wh{u} = [u;u';u''] \in \bR^{1+m+K}$,
which we later rename $u$,
and similarly for $v \in \bR$, $v' \in \bR^m$,
and $v'' \in \bR^K$,
the objective function simply reads $\eps \|\wh{u}\|_1 + \kappa \|\wh{v}\|_1$,
while the constraints read 
$\wh{u} + \wh{v} = [1;-a;0]$
and $[b,C,D]\wh{u}=0$.
The final transformation applied to arrive at the linear program of \eqref{LPforUB}
consists in introducing slack variables $r,s \in \bR^{1+m+K}$ such that $|\wh{u}| \le r$ and $|\wh{v}| \le s$. 
\epf

\brk
When the observation functionals and the quantity of interest are point evaluations,
a solution $a^{(\ul{t})}$ to the linear program \eqref{LPforUB} yields a linear functional $R^{(\ul{t})}: y \in \bR^m \mapsto \sum_{i=1}^m a^{(\ul{t})}_i y_i \in \bR$ which turns out to be a near-optimal recovery map.
Indeed,
remembering that the first step of the proof of Theorem~\ref{ThmGen} is valid for $a^{(\ul{t})}$,
which comes with atomic measures $\mu^{(\ul{t})}$ and $\nu^{(\ul{t})}$ such that $\mu^{(\ul{t})}+\nu^{(\ul{t})} = \rho - \sum_{i=1}^m a_i^{(\ul{t})} \la_i$
and 
$\int_{-1}^1 v(x) d\mu^{(\ul{t})}(x)=0$ for all $v \in \cP_n$, we derive that the error of $R^{(\ul{t})}$ over $\cK$ satisfies  
\be
 \sup_{f \in \cK} |Q(f) - R^{(\ul{t})}(L(f))|
\le \int_{-1}^1 \eps \, d|\mu^{(\ul{t})}| + \kappa \, d|\nu^{(\ul{t})}| 
= \beta^{(\ul{t})}.
\ee
This estimate matches the intrinsic error $\alpha^{\rm opt}$
with error at most $\beta^{(\ul{t})} - \alpha^{(N)}$.
This quantity, which is available after solving \eqref{OptTypeITrunc} and \eqref{LPforUB}, 
is small when the truncation parameter $N$ and the size of the grid $\ul{t}$ are large. 
\erk

\end{document}